\documentclass[12pt]{article}
\usepackage{amsmath,amssymb,amsfonts,graphicx,wrapfig,hyperref}
\usepackage[latin1]{inputenc}

\newtheorem{theorem}{Theorem}

\newtheorem{claim}{Claim}

\newcommand{\li}{{\mathrm{lim\,inf}}}
\newcommand{\ls}{{\mathrm{lim\,sup}}}

\title{On the sizes of large subgraphs of the binomial random graph}
\date{}
\author{J\'ozsef Balogh\footnote{Department of Mathematical Sciences, University of Illinois at Urbana-Champaign, Urbana, Illinois 61801, USA, 
and Moscow Institute of Physics and Technology, 9 Institutskiy per., Dolgoprodny, Moscow Region, 141701, Russian Federation.}
\and Maksim Zhukovskii\footnote{Moscow Institute of Physics and Technology, 9 Institutskiy per., Dolgoprodny, Moscow Region, 141701, Russian Federation; Adyghe State University, Caucasus mathematical center, ul. Pervomayskaya, 208, Maykop, Republic of Adygea, 385000, Russian Federation; The Russian Presidential Academy of National Economy and Public Administration, Prospect Vernadskogo, 84, bldg 2, Moscow, 119571, Russian Federation; Moscow Center for Fundamental and Applied Mathematics, Russian Federation.}
}

\begin{document}

\maketitle

\begin{abstract}
We consider the binomial random graph $G(n,p)$, where $p$ is a constant, and answer the following two questions. 

First, given $e(k)=p{k\choose 2}+O(k)$, what is the maximum $k$ such that a.a.s.~the binomial random graph $G(n,p)$ has an induced subgraph with $k$ vertices and $e(k)$ edges? We prove that this maximum is not concentrated in any finite set (in contrast to the case of a small $e(k)$). Moreover, for every constant $C>0$ and every $\omega_n\to\infty$, a.a.s.~the size of the concentration set belongs to $(C\sqrt{n/\ln n},\omega_n\sqrt{n/\ln n})$.

Second, given $k>\varepsilon n$, what is the maximum $\mu$ such that a.a.s.~the set of sizes of $k$-vertex subgraphs of $G(n,p)$ contains a full interval of length $\mu$? The answer is $\mu=\Theta\left(\sqrt{(n-k)n\ln{n\choose k}}\right)$.
\end{abstract}

\section{Introduction}

Let $\mathcal{F}$ be a family of graphs. Let $\mathcal{F}_k$ denote the set of graphs from $\mathcal{F}$ that have $k$ vertices. Let $X_n=X_n(\mathcal{F})$ be the maximum $k$ such that $G(n,p)$ has an induced subgraph isomorphic to some $F\in\mathcal{F}_k$. Below, we briefly discuss the main results on the asymptotic behaviour of $X_n$ (although we focus on constant $p$, we try to state all known results in the most general setting).\\

The first related result describes the asymptotic behaviour of the independence number (the maximum size of an independent set) and the clique number (the maximum size of a clique) of $G(n,p)$~\cite{BE_independent,M_independent_0,M_independent}. It states that, for arbitrary constant $p\in(0,1)$, there exists $f(n)$ such that asymptotically almost surely (a.a.s.) the clique number of $G(n,p)$ belongs to $\{f(n),f(n)+1\}$ (below, in such situations we say that the random variable {\it is 2-point concentrated}). By symmetry reasons, the same is true for the independence number. For the latter parameter, the same techniques work when $p=p(n)$ is large enough ($p\geq n^{-\varepsilon}$ for small enough constant $\varepsilon>0$), and for the clique number it works for small enough $p$ ($p\leq 1-n^{-\varepsilon}$). Certain improvements and generalizations of these results can be found in~\cite{Near_Concentration_Set,Rai}.\\

Clearly, the independence number and the clique number equal $X_n(\mathcal{F})$, where $\mathcal{F}$ is the family of empty graphs and the family of complete graphs respectively.

A natural question to ask is: what about other `common' graph sequences, such as paths, cycles, etc.? In~2018, Dutta and Subramanian \cite{Many_Induced_Graphs}, obtained 2-point concentration results for the family of simple paths and the family of simple cycles. Both results hold when $p\geq n^{-1/2}(\ln n)^2$.

Let us turn to larger graph families. The following families were considered by several researchers: trees, regular graphs, complete bipartite graphs and complete multipartite graphs~\cite{Palka_Trees,Bipartite, General_Induced}. Unfortunately, for all these families, it is still unknown, if 2-point concentration holds, or even if $m$-point concentration holds for some fixed number $m$. In 1983, Erd\H{o}s and Palka~\cite{Palka_Trees} proved that, for trees (i.e., $\mathcal{F}_k$ consists of all trees on $k$ vertices), $\frac{X_n}{\ln n}\stackrel{{\sf P}}\to \frac{2}{\ln[1/(1-p)]}$ as $n\to\infty$ (hereinafter, $\stackrel{{\sf P}}\to $ denotes the convergence in probability). In 1987, Ruci\'{n}ski~\cite{General_Induced} obtained a similar law of large numbers type general result for a quite wide class of graph families $\mathcal{F}_k$. In particular, from his result follows that: if $\mathcal{F}_k$ are sets of $ck(1+o(1))$-regular graphs, then $\frac{X_n}{\ln n}\stackrel{{\sf P}}\to\frac{2}{c\ln[1/p]+(1-c)\ln[1/(1-p)]}$ as $n\to\infty$. For several families of complete bipartite and multipartite graphs, similar results were obtained in~\cite{Bipartite, General_Induced}.\\

In 2012, Fountoulakis, Kang and McDiarmid \cite{Subgraphs_bounded_deg} considered families of graphs defined by constraints on the number of edges. More formally, given a sequence $e=e(k)$, $\mathcal{F}_k$ is the set of all graphs on $k$ vertices having {\it at most} $e(k)$ edges. The main result of \cite{Subgraphs_bounded_deg} states, in particular, the following. Let $n^{-1/3+\varepsilon}<p<1-\varepsilon$ for some $\varepsilon\in(0,1/3)$. Let $e=e(k)=o(\frac{pk\ln k}{\ln\ln k})$ be a sequence of non-negative integers. Then $X_n$ is 2-point concentrated.

It is easy to show, using the so-called second moment method, that a similar result holds for families of graphs having {\it exactly} $e$ edges: if $0\leq e(k)=O(k)$ (for sure, this bound can be improved, but there is no point to be very precise here), $\mathcal{F}_k$ is the set of all graphs on $k$ vertices with exactly $e(k)$ edges, then $X_n$ is 2-point concentrated. For the sake of convenience, let us denote this random variable $X_n$ by $\mathcal{X}_n(e)$.\\

One of the main goals of our study is to find a graph family $\mathcal{F}$ such that $X_n(\mathcal{F})$ is not concentrated within an interval of constant length (i.e.~there is no sequence $f(n)$ and fixed number $m$ such that a.a.s.~$X_n\in[f(n),f(n)+m]$). In particular, we want to find a sequence $e=e(k)$ such that $\mathcal{X}_n(e)$ is not concentrated within an interval of constant length. It is quite natural to check, if the `average' number of edges $e(k)=p{k\choose 2}+O(1)$ is appropriate (since the number of edges is integer, the right side should be integer as well --- this is why $O(1)$ appears). In other words, how many vertices should we remove from the random graph to make the number of edges equal to its expected number? Is this number of vertices tightly concentrated? We give the following answer for both questions for a much wider class of functions $e(k)$.

\begin{theorem} Let $p\in(0,1)$ be a constant. Let $e(k)={k\choose 2}p+O(k)$ be a sequence of non-negative integers.
\begin{enumerate}
\item[(i)] There exists a constant $t>0$ such that, for all constants $c>t$ and $C>2c+t$, we have
$$
0<\li_{n\to\infty}{\sf P}\left(n-C\sqrt{\frac{n}{\ln n}}<\mathcal{X}_n(e)<n-c\sqrt{\frac{n}{\ln n}}\right)\leq
$$
$$
 \ls_{n\to\infty}{\sf P}\left(n-C\sqrt{\frac{n}{\ln n}}<\mathcal{X}_n(e)<n-c\sqrt{\frac{n}{\ln n}}\right)<1.
$$
\item[(ii)] Let, for every sequence $m_k=O(\sqrt{k/\ln k})$ of non-negative integers, the following smoothness condition hold: $\left|\left(e(k)-{k\choose 2}p\right)-\left(e(k-m_k)-{k-m_k\choose 2}p\right)\right|=o(k)$. Then, for every constant $\varepsilon>0$, there exist constants $c_{\varepsilon},C_{\varepsilon}$ such that
$$
 \li_{n\to\infty}{\sf P}\left(n-C_{\varepsilon}\sqrt{\frac{n}{\ln n}}<\mathcal{X}_n(e)<n-c_{\varepsilon}\sqrt{\frac{n}{\ln n}}\right)>1-\varepsilon.
$$
\end{enumerate}
\label{th_smallf}

\end{theorem}

{\it Remark.} The first part of Theorem~\ref{th_smallf} implies that $\mathcal{X}_n(e)$ is not tightly concentrated. Moreover, under the additional smoothness condition in (ii), the size of the concentration set is $O_{{\sf P}}(\sqrt{\frac{\ln n}{n}})$ ($O_{{\sf P}}$ refers to the stochastic boundedness: for a sequence of random variables $\xi_n$ and a sequence of non-zero real numbers $a_n$, $\xi_n=O_{{\sf P}}(a_n)$ means that, for every $\varepsilon>0$, there exists $C>0$ such that ${\sf P}(|\xi_n/a_n|>C)<\varepsilon$), and this asymptotical bound is the best possible. The smoothness condition in (ii) holds for all $e(k)={k\choose 2}p+o(k)$.\\

This result is closely related to a study of possible sizes (i.e., number of edges) of subgraphs of the random graph that was initiated in 2009 by Alon and Kostochka~\cite{AlonKostochka}. 
Let us ask the following question. What is the maximum $\mu=\mu(k)$, $k\in\mathbb{N}$, such that a.a.s., for every $k$, the set of sizes of $k$-vertex induced subgraphs of $G(n,p)$ contains a full interval of length $\mu(k)$? In~\cite{AlonKostochka} it is proved that, for $k\leq 10^{-3}n$ and $p=1/2$, $\mu=\Omega(k^{3/2})$. 

This result is motivated by the following conjecture of Erd\H{o}s, Faudree and S\'{o}s (see~\cite{c-Ramsey_conj_1,c-Ramsey_conj_2}): {\it for every constant $c>0$, there exists a constant $b=b(c)>0$ so that if $G$ is a $c$-Ramsey graph on $n$ vertices, then the number of distinct pairs $(|V(H)|,|E(H)|)$, as $H$ ranges over all induced subgraphs of $G$, is at least $bn^{5/2}$} (an $n$-vertex graph is $c$-Ramsey if both its independence number and clique number are at most $c\ln n$; $V(H)$ and $E(H)$ denote the set of vertices and the set of edges of $H$ respectively). The result of \cite{AlonKostochka} immediately implies that the conjecture is true for almost all graphs. In 2019, the conjecture was proved by Kwan and Sudakov~\cite{KwanSudakov}.

Extending results of~\cite{AlonKostochka}, we 
get asymptotically close upper and lower bounds (that differ by a multiplicative constant factor) on $\mu$ for $k>\varepsilon n$.

\begin{theorem} Let $p\in(0,1)$ be a constant. Let $\varepsilon>0$ be an arbitrary small constant.
\begin{enumerate}
\item[(i)] There exists a constant $q>0$ such that a.a.s., for every $k\in\{\lfloor\varepsilon n\rfloor,\ldots,n-1\}$, the set of sizes of induced $k$-vertex subgraphs of $G(n,p)$ contains a full interval of length at least $qk\sqrt{\frac{n-k}{n}\ln{n\choose k}}$. Moreover, a.a.s., for every $k\in\{1,\ldots,\lfloor\varepsilon n\rfloor-1\}$, the set of sizes of induced $k$-vertex subgraphs of $G(n,p)$ contains a full interval of length at least $qk^{3/2}$. 

\item[(ii)] There exists a constant $Q>0$ such that a.a.s., for every $k\in\{1, \ldots,n-1\}$, the set of sizes of induced $k$-vertex subgraphs of $G(n,p)$ does not contain any full interval of length at least $Q k\sqrt{\frac{n-k}{n}\ln{n\choose k}}$. 
\end{enumerate}
\label{AK}
\end{theorem}

Therefore, for $k\geq\varepsilon n$, $\mu=\Theta\left(\sqrt{(n-k)n\ln{n\choose k}}\right)$. For $k<\varepsilon n$, $\mu\in\left[qk^{3/2},Qk^{3/2}\sqrt{\ln(n/k)}\right]$ for some constants $q,Q$. The latter lower bound ($\mu\geq qk^{3/2}$) follows immediately from the result of~\cite{AlonKostochka} since their proof works for arbitrary constant $p$. Notice that, for $k=n-o(n)$, the value of $k\sqrt{\frac{n-k}{n}\ln{n\choose k}}$ becomes much smaller than $k^{3/2}$. For such $k$, the techniques of Alon and Kostochka~\cite{AlonKostochka} does not work since one should switch from analysis of spanned edges of induced $k$-subgraphs to analysis of edges having at least one end-point outside.

Notice that, for all $k<\frac{2-\delta}{\max\{\ln(1/p),\ln(1/(1-p))\}}\ln n$, the exact value of $\mu(k)$ is known: $\mu(k)={k\choose 2}+1$ since a.a.s., for every such $k$ and every graph $F$ on $k$ vertices, there is an induced subgraph in $G(n,p)$ isomorphic to $F$ (this is a simple exercise that can be solved using the second moment method; for $p=1/2$, it appears as exercise 1 in~\cite{AlonSpencer}). In a similar manner, one may prove that for every $k=\Theta(\ln n)$ our upper bound has right order of magnitude: $\mu(k)=\Theta(k^2)=\Theta\left(k\sqrt{\frac{n-k}{n}\ln{n\choose k}}\right)$.


\section{Preliminaries}
\label{sec:pre}

\subsection{Concentration of binomial random variables}

Let $\Phi(x)=\int_{-\infty}^x \frac{1}{\sqrt{2\pi}}e^{-t^2/2}dt$. Consider a binomial random variable $\xi$ with parameters $N$ and $p$. Then, by the DeMoivre--Laplace theorem (see, e.g.,~\cite{Bol_Degree2}~and~\cite{Feller}), for $h=o(N^{1/6})$, 
\begin{equation}
{\sf P}\left[\xi\leq Np+h\sqrt{Np(1-p)}\right]\sim\Phi(h),
\label{moivre_laplace}
\end{equation} 
\begin{equation}
\text{for integer }Np+h\sqrt{Np(1-p)},\quad
{\sf P}\left[\xi= Np+h\sqrt{Np(1-p)}\right]\sim\frac{1}{\sqrt{2\pi Np(1-p)}}e^{-h^2/2}
\label{moivre_laplace_equal}
\end{equation}
as $N\to\infty$.
In our proofs, we multiple times use the following relation: 
\begin{equation}
1-\Phi(x)\sim\frac{1}{\sqrt{2\pi}x}e^{-x^2/2}\text{ as }x\to\infty
\label{phi_approximation}
\end{equation} 
(see relation ($1^{\prime}$) in \cite{Bol_Degree2}).

We also need the following version of the Chernoff bound (see, e.g., \cite[Theorem 2.1]{Janson}: for $t\geq 0$,
$$
 {\sf P}(\xi\geq Np+t)\leq\exp\left(-\frac{t^2}{2(Np+t/3)}\right),\quad
 {\sf P}(\xi\leq Np-t)\leq\exp\left(-\frac{t^2}{2Np}\right).
$$

\subsection{Graph notations}

Given a graph $\Gamma$ and a set $U\subset V(\Gamma)$, we denote by $\Gamma|_U$ the subgraph of $\Gamma$ induced on the set $U$. We call the number of edges of $\Gamma$ having vertices in $U$ {\it the degree of $U$} and denote it $\delta_{\Gamma}(U)$ (i.e., $\delta(U)=|\{\{u,v\}\in E(\Gamma):\,\text{either }u\in U,\text{ or }v\in U\}|$) or simply $\delta(U)$ for $\Gamma=G(n,p)$.

We also use notations $v(\Gamma)$ and $e(\Gamma)$ for the number of vertices and the number of edges in $\Gamma$ respectively; $\Delta[\Gamma]$ denotes the maximum degree of $\Gamma$.

\subsection{Degrees in random graphs}

As usual, the vertex set of $G(n,p)$ is $[n]:=\{1,\ldots,n\}$. We will use the following fact: a.a.s.~the maximum degree of $G(n,p)$ is at most $pn+\sqrt{2p(1-p)n\ln n}$~\cite{Bol_Degree}.

Let us now state a very helpful claim on the number of vertices in $G(n,p)$ having large degrees. Consider sequences 
$0<\alpha=\alpha(n)\leq 2\ln n-5\ln\ln n$ and $\tilde\alpha=\tilde\alpha(n)$ such that $|\tilde\alpha|\leq\sqrt{\ln n-\frac{5}{2}\ln\ln n}$.
Set $\zeta(\alpha)=np+\sqrt{\alpha np(1-p)}$, $\tilde\zeta(\tilde\alpha)=\lfloor (n-1)p+\tilde\alpha\sqrt{np(1-p)}\rfloor$. Let $\varepsilon>0$. Denote 
$$
g_{\varepsilon}:=g_{\varepsilon}(n)=(1-\varepsilon)\frac{n}{\sqrt{2\alpha\pi}}e^{-\frac{\alpha}{2}},\quad
\tilde g_{\varepsilon}:=\tilde g_{\varepsilon}(n)=(1-\varepsilon)\frac{\sqrt{n}}{\sqrt{2\pi p(1-p)}}e^{-\frac{\tilde\alpha^2}{2}}.
$$

\begin{claim}
In $G(n,p)$, with probability $o(1/n^2)$,
\begin{enumerate}
\item there are at most $g_{\varepsilon}$ vertices having degrees greater than $\zeta(\alpha)$,
\item there are at most $\tilde g_{\varepsilon}$ vertices having degrees equal to $\tilde\zeta(\tilde\alpha)$.
\end{enumerate}
\label{claim:general_degrees}
\end{claim}

{\it Proof.} Let $Y$ and $Z$ be the number of vertices in $G(n,p)$ having degrees greater than $\zeta(\alpha)$ and equal to $\tilde\zeta(\tilde\alpha)$ respectively. Consider a normal random variable $\mathbf{p}$ with mean $p$ and variance $\frac{p(1-p)}{n(n-1)}$, truncated to $(0,1)$. Assume that, for every $x\in(0,1)$, independent binomial random variables $X_1(x),\ldots,X_n(x)\sim\mathrm{Bin}(n-1,x)$ are defined independently of $\mathbf{p}$. Assume also  that all the above random variables are defined on the same probability space as $G(n,p)$. Let $d_1,\ldots,d_n$ be the degrees of the vertices $1,\ldots,n$ in $G(n,p)$. Then (see, e.g.,~\cite[Proposition 1.2]{Wormald}), for every set $S\subset\mathbb{R}^n$, 
$$
{\sf P}((d_1,\ldots,d_n)\in S)\leq {\sf P}\biggl[(X_1(\mathbf{p}),\ldots,X_n(\mathbf{p}))\in S\left|\sum_{i=1}^n X_i(\mathbf{p})\text{ is even}\right.\biggr]+o(1/n^2).
$$
Let $f$ be the density of $\mathbf{p}$. Then, from (\ref{phi_approximation}), 
$$
{\sf P}\left(|\mathbf{p}-p|>\frac{\ln n}{e n}\right)=\int_0^{p-\frac{\ln n}{n}}f(x)dx+\int_{p+\frac{\ln n}{n}}^1f(x)dx=
$$
$$
\frac{1}{\Phi(1)-\Phi(0)}\left[\int_{-\sqrt{\frac{p}{1-p}n(n-1)}}^{-\ln n\sqrt{\frac{1-1/n}{p(1-p)}}}+\int_{\ln n\sqrt{\frac{1-1/n}{p(1-p)}}}^{\sqrt{\frac{1-p}{p}n(n-1)}}\right]\frac{1}{\sqrt{2\pi}}e^{-x^2/2}dx=o\left(\frac{1}{n^2}\right).
$$
Also, since $\sum_{i=1}^{n}X_i(x)$ has binomial distribution with parameters $n(n-1)$ and $x$, we get that,  for every $x\in(0,1)$,
$$
 {\sf P}\left(\sum_{i=1}^{n}X_i(x)\text{ is even}\right)=\frac{1}{2}+\frac{1}{2}(1-2x)^{n(n-1)}.
$$
Let $Y(x)$ ($Z(x)$) be the number of $i\in[n]$ such that $X_i(x)>\zeta(\alpha)$ ($X_i(x)=\tilde\zeta(\tilde\alpha)$ resp.). Then, from the above,
$$
 {\sf P}\left(Y\leq g_{\varepsilon}\right)\leq \int_0^1{\sf P}\left(Y(x)\leq g_{\varepsilon}\left|\sum_{i=1}^{n}X_i(x)\text{ is even}\right.\right)f(x)dx+o\left(\frac{1}{n^2}\right)\leq
$$
\begin{equation} 
 \frac{{\sf P}\left(Y(p-\ln n/ n)\leq cn\right)}{\frac{1}{2}+\frac{1}{2}(1-2p-2\ln n/ n)^{n(n-1)}}+o\left(\frac{1}{n^2}\right)=2{\sf P}\left(Y(p-\ln n/ n)\leq g_{\varepsilon}\right)+o\left(\frac{1}{n^2}\right).
\label{from_random_p_to_non_random}
\end{equation}
By~(\ref{moivre_laplace})~and~(\ref{phi_approximation}), for every $i\in[n]$,
$$ 
{\sf P}(X_i(p-\ln n/ n)>\zeta(\alpha))
=\frac{1}{\sqrt{2\alpha\pi}}e^{-\frac{\alpha}{2}}(1+o(1)).
$$
Finally, from the Chernoff bound, we get
$$ {\sf P}\left(Y(p-\ln n/ n)\leq g_{\varepsilon}\right) \leq
 \exp\left[-\varepsilon^2\frac{n}{2\sqrt{2\alpha\pi}}e^{-\frac{\alpha}{2}}\right]
 =\exp\left[-\Omega(\ln^2 n)\right]=o(1/n^2).
$$
In the same way, since for every $x\in[p-\ln n/ n,p+\ln n/ n]$, by~(\ref{moivre_laplace_equal}),
$$ 
{\sf P}(X_i(x)=\tilde\zeta(\tilde\alpha))
=\frac{1}{\sqrt{2\pi np(1-p)}}e^{-\frac{\tilde\alpha^2}{2}}(1+o(1))
$$
we get that, for every such $x$,
$$
{\sf P}\left(Z(x)\leq \tilde g_{\varepsilon}\right)\leq
\exp\left[-\varepsilon^2\frac{\sqrt{n}}{2\sqrt{2\pi p(1-p)}}e^{-\frac{\tilde\alpha^2}{2}}\right](1+o(1))=\exp\left[-\Omega(\ln^{5/4} n)\right]=o\left(\frac{1}{n^2}\right).
$$
Therefore,
$$
{\sf P}\left(Z\leq \tilde g_{\varepsilon}\right) \leq 2{\sf P}\left(Z(p)\leq \tilde g_{\varepsilon}\right)(1+o(1))+o\left(\frac{1}{n^2}\right)=o\left(\frac{1}{n^2}\right).\quad\Box
$$



\section{Proof of Theorem~\ref{th_smallf}}


Denote $f(k)=e(k)-{k\choose 2}p$. Below in the proof, we assume that $Q\in\mathbb{R}$ is such that $-Qk\leq f(k)\leq Qk$ for all $k\in\mathbb{N}$. \\

The proof is divided into three parts. 

In Section~\ref{sec:intervals}, we consider several bounds on the number of edges in $G(n,p)$ that hold with positive asymptotic probabilities. That is, we consider two intervals $I_n^1,I_n^2$ and a set $I_n(\varepsilon)$ such that the left bound $p{n\choose 2}+(a_2+Q)n$ of $I_n^2$ is bigger than the right bound $p{n\choose 2}+(b_1-Q)n$ of $I_n^1$. 
All intervals are of sizes $O(n)$, the asymptotical probability that the number of edges is inside $I_n^j$, $j\in\{1,2\}$, is positive, and the probability of the same event but for $I_n(\varepsilon)$ is bigger than $1-\varepsilon$. 

In Section~\ref{sec:upper}, we obtain upper bounds on $\mathcal{X}_n(e)$. First, we assume that $e(G(n,p))\in I_n^i$ and obtain upper bounds $B_i=n-c_i\sqrt{n/\ln n}$. Second, we assume that $e(G(n,p))\in I_n(\varepsilon)$, and obtain an upper bound $B(\varepsilon)=n-c\sqrt{n/\ln n}$. 

In Section~\ref{sec:lower}, we obtain lower bounds on $\mathcal{X}_n(e)$. First, we assume that $e(G(n,p))\in I_n^i$ and obtain lower bounds $A_i=n-C_i\sqrt{n/\ln n}$. Second, we assume that $e(G(n,p))\in I_n(\varepsilon)$ and obtain a lower bound $A(\varepsilon)=n-C\sqrt{n\ln n}$. 

Combining the second and the third part, we obtain that, first, the lower bound $A_2$ is bigger than the upper bound $B_1$ whenever $a_2>2b_1$. This finishes the proof of Theorem~\ref{th_smallf}.(i). Second, since both bounds $A(\varepsilon)$ and $B(\varepsilon)$ are true with asymptotical probabilities at least $1-\varepsilon$, we get Theorem~\ref{th_smallf}.(ii).

\subsection{Bounds on the number of edges}
\label{sec:intervals}

 Fix real numbers $a_1<b_1<a_2<b_2$ such that $a_1>0$, $b_1>a_1+2Q$, $a_2>2b_1$, $b_2>a_2+2Q$. Consider the sets 
$$
I_n^1=\left(p{n\choose 2}+(a_1+Q)n,\, p{n\choose 2}+(b_1-Q)n\right),
$$
$$
I_n^2=\left(p{n\choose 2}+(a_2+Q)n,\, p{n\choose 2}+(b_2-Q)n\right).
$$
Let $\gamma>0$ be such that, for $n$ large enough, 
\begin{equation}
 \min\left\{{\sf P}(e(G(n,p))\in I_n^1),\,{\sf P}(e(G(n,p))\in I_n^2)\right\}>\gamma.
\label{two_intervals}
\end{equation}
Such $\gamma$ exists since $e(G(n,p))\sim$Bin$({n\choose 2},p)$ and $p$ is constant, see Section~\ref{sec:pre}.

Moreover, for every $\varepsilon>0$, choose $a=a(\varepsilon)$ and $b=b(\varepsilon)$ such that, for $n$ large enough,
$$
 {\sf P}\left(e(G(n,p))\in I_n(\varepsilon)\right)>1-\varepsilon, \text{ where } 
$$
\begin{equation} 
 I_n(\varepsilon)=\left(p{n\choose 2}-(b-Q)n,p{n\choose 2}+(b-Q)n\right)\setminus 
 \left[e(n)-an,e(n)+an\right].
\label{third_interval}
\end{equation}

\subsection{Upper bounds on $\mathcal{X}_n(e)$}
\label{sec:upper}

Consider a sequence of integers $m=m(n)\leq\frac{c}{\sqrt{2p(1-p)}}\sqrt{\frac{n}{\ln n}}$. Denote $M=M(m)={m\choose 2}+m(n-m)$ the maximum possible degree of an $m$-set. Then, for a fixed $m$-set, the expected value of its degree equals $pM$. Consider the random variable 
\begin{equation}
Y_m=\max_{U\in{[n]\choose m}}\delta(U),
\label{def_Y}
\end{equation}
where ${[n]\choose m}$ is, as usual, the set of $m$-element subsets of $[n]$.

Then $Y_1=\Delta[G(n,p)]$ is the maximum degree of $G(n,p)$.  Since $Y_1<pn+\sqrt{2p(1-p)n\ln n}$ holds a.a.s. (see Section~\ref{sec:pre}), we immediately get that, a.a.s. 
$$
Y_m\leq mY_1<mpn+m\sqrt{2p(1-p)n\ln n}=Mp+m\sqrt{2p(1-p)n\ln n}+o(n),
$$
that is, 
a.a.s.~$Y_m<Mp+c n+o(n)$. Under the assumption that $e(G(n,p))>p{n\choose 2}+(a_i+Q)n$, we should remove a set of $m$ vertices with degree at least $a_i n$ bigger than the average value $pM$ in order to obtain a graph with at most $Qn$ edges more than the average value $p{n-m\choose 2}$.  
Thus, if $c<a_i$, a.a.s.~we cannot reach the desired number of edges by removing an $m$-set. Therefore, for every $\delta>0$, from~(\ref{two_intervals}), we get that 
\begin{equation}
{\sf P}\left(\mathcal{X}_n(e)<n-\frac{a_i(1-\delta)}{\sqrt{2p(1-p)}}\sqrt{\frac{n}{\ln n}}\right)>\gamma-o(1)
\label{upper}
\end{equation} 
for $i\in\{1,2\}$.

If $|e(G(n,p))-e(n)|>a(n)$, then, in order to obtain exactly $e(n-m)$ edges, we should remove a set of $m$ vertices with degree at least $e(n)-e(n-m)+a(n)=pM+f(n)-f(n-m)+a(n)$. Since $|f(n)-f(n-m)|=o(n)$, in the same way, from~(\ref{third_interval}), we get that 
\begin{equation}
{\sf P}\left(\mathcal{X}_n(e)<n-\frac{a(1-\delta)}{\sqrt{2p(1-p)}}\sqrt{\frac{n}{\ln n}}\right)>1-\varepsilon-o(1).
\label{upper_second}
\end{equation}

\subsection{Lower bounds on $\mathcal{X}_n(e)$} 
\label{sec:lower}

This part of the proof is divided into five parts. The overall idea is to use a small set of vertices (we extract it in Section~\ref{sec:parts}) to make the number of edges precisely $e(k)$. This small set appears helpful after the major part of extra edges is destroyed. More precisely, having $(b+Q)n$ edges more than the average, we can easily destroy extra $bn$ edges by removing a set of $O(\sqrt{n/\ln n})$ vertices. We do that in Section~\ref{sec:estimate}. But this is far from what we need since $f$ may differ a lot from its bound $Q$. In Section~\ref{sec:major}, we show how to reduce the number of extra edges up to $O(\sqrt{n\ln n})$. We use the supplementary small set in Sections~\ref{sec:first_part} and~\ref{sec:second_part} where we get the precise number of edges in two steps exploiting two equal parts of the set.

\subsubsection{Extracting a supplementary part} 
\label{sec:parts}

Let $n_0=\left\lfloor\frac{\sqrt{n}}{\ln n}\right\rfloor$, $\tilde n=n-2n_0$. Consider the partition $[n]=\{1,\ldots,2n_0\}\sqcup \tilde V_{\tilde n}$, where $\tilde V_{\tilde n}=\{2n_0+1,\ldots,n\}$. Divide the {\it supplementary set} $\{1,\ldots,2n_0\}$ into two disjoint parts of equal sizes $V_1=\{1,\ldots,n_0\}$ and $V_2=\{n_0+1,\ldots,2n_0\}$. Denote by $G_{\tilde n}$ the subgraph of $G(n,p)$ induced by $\tilde V_{\tilde n}$.


\subsubsection{Estimating from above the number of vertices we need to remove} 
\label{sec:estimate}

Fix $c>0$. By Claim~\ref{claim:general_degrees},  
a.a.s.~there are more than $m_c:=c\sqrt{\frac{n}{\pi\ln n}}$ vertices having degrees bigger than $np+\sqrt{\frac{1}{2}np(1-p)\ln n}+O(\sqrt{n})$. So, a.a.s. 
$$
Y_{m_c}[G_{\tilde n}]>\left(np+\sqrt{\frac{1}{2}np(1-p)\ln n}+O(\sqrt{n})\right)m_c-(m_c)^2=
M(m_c)p+cn\sqrt{\frac{p(1-p)}{2\pi}}+o(n),
$$
where the random variable $Y_{m_c}$ is defined in~(\ref{def_Y}), and $Y_{m_c}[G_{\tilde n}]$ is defined on the random graph $G_{\tilde n}$ (i.e. the maximum value of $\delta_{G_{\tilde n}}(U)$ over $U\in{{\tilde V_{\tilde n}}\choose m}$).

From the above, it follows that, in order to remove extra $Cn$ (or more) edges, we need to remove $m\leq C\sqrt{\frac{2n}{p(1-p)\ln n}}$ vertices. In the next section, we prove that this upper bound for the number of deleted vertices $m$ also works if we want to get the number of edges very close to the desired value $e(n-m)$ (the error term is $O(\sqrt{n\ln n})$).

\subsubsection{Removing a major part of extra edges}
\label{sec:major} 

For $m\in\mathbb{N}$ and $U\subset\tilde V_{\tilde n}$, set 
$$
\tilde M(m)=m(\tilde n-m)+{m\choose 2},\quad \overline{U}=\tilde V_{\tilde n}\setminus U,\quad \tilde\delta(U)=\delta_{G_{\tilde n}}(U).
$$
Moreover, let $E_{\tilde n}:=e(G_{\tilde n})-{\tilde n\choose 2}p$. From~(\ref{two_intervals}), ${\sf P}(E_{\tilde n}\in((a_i+Q)\tilde n,\,(b_i-Q)\tilde n))>\gamma$ for $n$ large enough and $i\in\{1,2\}$. 

Let us describe an algorithm of constructing a set of $m=O(\sqrt{n/\ln n})$ vertices $U\subset\tilde V_{\tilde n}$ such that $G_{\tilde n}|_{\overline{U}}$ has ${{\tilde n}-m\choose 2}p+f(\tilde n-m)+O(\sqrt{n\ln n})$ edges. 


At step $1$, $U_1=\{v_1\}$ where $v_1$ has maximum degree in $G_{\tilde n}$. If 
$$
\tilde\delta(U_1)>p\tilde M(1)+E_{\tilde n}-f(\tilde n-1),
$$ 
then the algorithm terminates, and $U=U_0:=\varnothing.$ 

Assume that, at step $i\geq 1$, we have a set $U_i$ of $i$ vertices. If the algorithm still runs, then consider the set $U_{i+1}=U_i\cup\{v_{i+1}\}$ of $i+1$ vertices having maximum degrees in $G_{\tilde n}$. If 
$$
\tilde\delta(U_{i+1})>p\tilde M(i+1)+E_{\tilde n}-f(\tilde n-i-1),
$$
then the algorithm terminates, and $U=U_i$. 

By results from Section~\ref{sec:estimate}, with probability at least $\gamma-o(1)$, the algorithm terminates after $O(\sqrt{n/\ln n})$ steps.\\

Let us prove that the algorithm gives a set of vertices $\overline{U}$ inducing a graph with the desired amount of edges but $O(\sqrt{n\ln n})$.

Let $i=O(\sqrt{\frac{n}{\ln n}})$. Let us estimate from above 
$$
(\tilde\delta(U_{i+1})-p\tilde M(i+1))-(\tilde\delta(U_i)-p \tilde M(i))=\tilde\delta(U_{i+1})-\tilde\delta(U_i)-p(\tilde n-i-1).
$$
Obviously, it cannot be bigger than $\Delta[G_{\tilde n}]-p(\tilde n-i-1)$. But the latter is larger than $2\sqrt{np(1-p)\ln n}+O(\sqrt{n})$ with probability $o(\frac{1}{n})$. Indeed, by~(\ref{moivre_laplace}),~(\ref{phi_approximation}) and the union bound,
$$
 {\sf P}(\Delta[G_{\tilde n}]>\tilde np+2\sqrt{\tilde np(1-p)\ln \tilde n})\leq n(1+o(1))\int_{2\sqrt{\ln\tilde n}}^{\infty}\frac{1}{\sqrt{2\pi}}e^{-x^2/2}dx\sim\frac{1}{2n\sqrt{2\pi\ln n}}.
$$

By using the described algorithm we conclude that here exists $\varepsilon_n\to 0$ such that, for every $\delta>0$ and $n$ large enough,
\begin{align*}
 {\sf P}\biggl[\exists U\subset\tilde V_{\tilde n}&\,\exists \alpha_0=\alpha_0(n)\in(3,6): \\
 &|U|=m\leq (b_i+\delta)\sqrt{\frac{2}{p(1-p)}}\sqrt{\frac{n}{\ln n}},\\
& \left|E\left(G_{\tilde n}|_{\overline{U}}\right)\right|={\tilde n-m\choose 2}p+f(\tilde n-m)+\alpha_0(n)\sqrt{np(1-p)\ln n}\biggr]\geq\gamma-\varepsilon_n.
\end{align*}


\subsubsection{Exploiting the first part of the supplementary set}  
\label{sec:first_part}

Here, we assume that the above algorithm constructs the desired set $U$ (this happens with probability at least $\gamma-o(1)$) of size $m$, and all the events below are conditioned by this event. Since all the events below are defined by edges chosen independently of $G(n,p)|_{\tilde V_{\tilde n}}$, we are still working with independent Bernoulli random variables.\\

For $V\subset V_1$, denote $\delta_{\overline{U}}(V)=\sum_{v\in V}\delta_{\overline{U}}(v)$ where $\delta_{\overline{U}}(v)$ is the  number of neighbors of $v$ in $\overline{U}$. From Claim~\ref{claim:algorithm} stated below ($V=V_1$, $W=\overline{U}$, $\ell=3$, $\beta=0.49$) it follows that a.a.s. there exists a set $V_1^0\subset V_1$ of $h=7$ vertices such that its recovering corrects the deviation from the desired number of edges up to $o(\sqrt{n})$, i.e. the difference between $\left|E\left(G(n,p)|_{\overline{U}\sqcup V_1^0}\right)\right|$ and the desired value $e(|\overline{U}\sqcup V_1^0|)$ equals $o(\sqrt{n})$. Define $\tilde U=\overline{U}\sqcup V_1^0$.

\begin{claim}
Let $\ell\in\mathbb{N}$, $\varepsilon>0$, $\beta>0$. Let $h=\left\lfloor\frac{\ell+3}{\sqrt{2\beta}}\right\rfloor+1$. Let $V\subset[n]$, $W\subset[n]\setminus V$ be such that $|V|\geq n^{\beta}$ and $|W|\geq\varepsilon n$. Then 
\begin{multline*}
 {\sf P}\biggl[\forall\alpha_0\in(3,3+\ell)\,\,\exists v_1,\ldots,v_h\in V\\
 \delta_W(\{v_1,\ldots,v_h\})=h|W|p-\alpha_0\sqrt{p(1-p)|W|\ln|W|}+O\left(\frac{\sqrt{n}}{\sqrt[4]{\ln n}}\right)\biggr]=1-o\left(\frac{1}{n}\right).
\end{multline*}
\label{claim:algorithm}
\end{claim}

The claim is proven in Section~\ref{claim:algorithm_proof}.



\subsubsection{Exploiting the second part of the supplementary set} 
\label{sec:second_part}

Here, we exploit the set $V_2$ and finish the construction of an induced graph with $k$ vertices and exactly $e(k)$ edges. The existence of the construction follows from the claim stated below (the proof of the claim is given in Section~\ref{claim:two_sets_proof}).

\begin{claim}
Let $\varphi:=\varphi(n)=o(\sqrt{n})$ be a sequence of real number, $\ell\in\mathbb{N}$, $\varepsilon>0$. Let $V\subset[n]$, $U\subset[n]\setminus V$ be such that $|V|\geq n^{\varepsilon}$, $|V|^{\ell}\geq n^{1/2+\varepsilon}$ and $|U|\geq\varepsilon n$. Then there are no vertices $w_1,\ldots,w_{\ell}$ in $V$ such that the number of edges in $G(n,p)|_{U\cup\{w_1,\ldots,w_{\ell}\}}$ having at least one end-point in $\{w_1,\ldots,w_{\ell}\}$ equals $\left\lfloor\left({\ell\choose 2}+\ell|U|\right)p+\varphi\right\rfloor$ with probability $o(1/n^2)$.
\label{claim:two_sets}
\end{claim}

Let us apply Claim~\ref{claim:two_sets} to $V=V_2$, $U=\tilde U$, $\ell=2$ and $\varphi=|E(G(n,p)|_{\tilde U})|-e(|\tilde U|)$ (notice that $\varphi$ is random, so we also need to apply the union bound over its possible values).

Since the difference of number of edges in $G(n,p)|_{\tilde U}$ from $e(|\tilde U|)$ is $o(\sqrt{n})$, then by Claim~\ref{claim:two_sets} and the union bound (over all possible values of the difference between the number of edges in $G(n,p)|_{\tilde U}$ and $e(|\tilde U|)$), a.a.s.~there exist two vertices $w_1,w_2$ in $V_2$ such that the graph $G(n,p)|_{\tilde U\cup\{w_1,w_2\}}$ has 
$$
k:=\tilde n-m+h+2\geq n-(b_i+\delta)\sqrt{\frac{2}{p(1-p)}}\sqrt{\frac{n}{\ln n}}+O\left(\frac{\sqrt{n}}{\ln n}\right)
$$ 
vertices and exactly $e_k$ edges.\\

Finally, we get that, for every $\delta>0$,
\begin{equation}
\li_{n\to\infty}{\sf P}\left(\mathcal{X}_n(e)>n-\frac{b_i(1+\delta)\sqrt{2}}{\sqrt{p(1-p)}}\sqrt{\frac{n}{\ln n}}\right)>\gamma
\label{lower}
\end{equation}
for both $i\in\{1,2\}$.

First, let $i=1$. Both (\ref{upper}) and (\ref{lower}) are obtained from (\ref{two_intervals}) (i.e., both events are intersections of one common event having probability bigger than $\gamma$ with events that hold a.a.s.). Therefore,
$$
 \li_{n\to\infty}{\sf P}\left(n-\frac{b_1(1+\delta)\sqrt{2}}{\sqrt{p(1-p)}}\sqrt{\frac{n}{\ln n}}<\mathcal{X}_n(e)<n-\frac{a_1(1-\delta)}{\sqrt{2p(1-p)}}\sqrt{\frac{n}{\ln n}}\right)>\gamma.
$$
Second, let $i=2$. Since $a_2>2b_1$, from (\ref{upper}) and (\ref{lower}), we get that
$$
 \ls_{n\to\infty}{\sf P}\left(n-\frac{a_2(1-\delta)}{\sqrt{2p(1-p)}}\sqrt{\frac{n}{\ln n}}\leq \mathcal{X}_n(e)\leq n-\frac{b_1(1+\delta)\sqrt{2}}{\sqrt{p(1-p)}}\sqrt{\frac{n}{\ln n}}\right)<1-2\gamma.
$$
Putting $t=\frac{2Q\sqrt{2}}{\sqrt{p(1-p)}}$, we finish the proof of Theorem~\ref{th_smallf}.(i).\\

In the same way, from (\ref{third_interval}), we get that
$$
\li_{n\to\infty}{\sf P}\left(\mathcal{X}_n(e)>n-\frac{b(1+\delta)\sqrt{2}}{\sqrt{p(1-p)}}\sqrt{\frac{n}{\ln n}}\right)>1-\varepsilon.
$$
Together with (\ref{upper_second}), this finishes the proof of Theorem~\ref{th_smallf}.(ii).\\

\subsubsection{Proof of Claim~\ref{claim:algorithm}}
\label{claim:algorithm_proof}

For a subset $V_0\subseteq V$, let 

\begin{itemize}
\item $p_h(V_0)$ be the probability that all but at most $h-1$ vertices of $V_0$ have more than $|W|p-\frac{3+\ell}{h}\sqrt{p(1-p)|W|\ln|W|}$ neighbors in $W$,

\item $\tilde p_h(V_0)$ be the probability that all vertices of $V_0$ have more than $|W|p-\frac{3}{h}\sqrt{p(1-p)|W|\ln|W|}$ neighbors in $W$.
\end{itemize}

Set $\kappa:=|V_0|$. By~(\ref{moivre_laplace}) and~(\ref{phi_approximation}), 
$$
p_h(V_0)=\sum_{i=0}^{h-1}{\kappa\choose i}\left(n^{-\frac{(3+\ell)^2}{2h^2}+o(1)}\right)^i\left(1-n^{-\frac{(3+\ell)^2}{2h^2}+o(1)}\right)^{\kappa-i},\quad
\tilde p_h(V_0)=\left(1-n^{-\frac{9}{2h^2}+o(1)}\right)^{\kappa}.
$$
Since $|V|\geq n^{\beta}$, $p_h(V)$ approaches $0$. Let $0<\tilde\beta<\frac{9}{2h^2}$. Set $V_0:=\{1,\ldots,\lfloor n^{\tilde\beta}\rfloor\}$. Then $\tilde p_h(V_0)\to 1$.\\

Let us describe an algorithm of constructing a sequence of $h$-sets of vertices from $V$ such that, for every $\alpha_0$, at least one of the sets in this sequence is the desired one.\\

Start from $V_0^1=V_0$. Consider $h$ vertices $u_1^1,\ldots,u_h^1$ of $V_0^1$ that have minimum number of neighbors in $W$. If $\delta_{W}(\{u_1^1,\ldots,u_h^1\})\leq h|W|p-(3+\ell)\sqrt{|W|p(1-p)\ln|W|}$, then the algorithm terminates. Clearly, the probability that the algorithms terminates at the first step is at most $1-\tilde p_h(V_0)\to 0$ as $n\to\infty$.

Let, at step $\kappa\geq 1$, the set $V_0^{\kappa}$ be constructed and the algorithm still works. Then, at step $\kappa+1$, consider $V_0^{\kappa+1}=V_0^{\kappa}\cup \{\lfloor n^{\tilde\beta}\rfloor+\kappa-1\}$ and choose $h$ vertices $u_1^{\kappa+1},\ldots,u_h^{\kappa+1}$ from it  that have minimum number of neighbors in $W$. If $\delta_{W}(\{u_1^{\kappa+1},\ldots,u_h^{\kappa+1}\})\leq h|W|p-(3+\ell)\sqrt{|W|p(1-p)\ln|W|}$, then the algorithm terminates.

It remains to prove that a.a.s., for every $\kappa$, $\delta_{W}(\mathbf{u}_{\kappa-1})-\delta_{W}(\mathbf{u}_{\kappa})=O\left(\frac{\sqrt{n}}{\sqrt[4]{\ln n}}\right)$ where $\mathbf{u}_{\kappa}=\{u_1^{\kappa},\ldots,u_h^{\kappa}\}$.\\

Clearly, the sets $\mathbf{u}_{\kappa}$ and $\mathbf{u}_{\kappa-1}$ have at least $h-1$ vertices in the intersection. Let $u^{\kappa}_1,\ldots,u^{\kappa}_{h+1}$ be the vertices of $V_0^{\kappa}$ having minimum number of neighbors in $W$, and $\delta_{W}(u^{\kappa}_1)\leq\ldots\leq \delta_{W}(u^{\kappa}_{h+1})$. Then, 
\begin{equation}
\delta_{W}(\mathbf{u}_{\kappa-1})-\delta_{W}(\mathbf{u}_{\kappa})\leq \delta_{W}(u^{\kappa}_{h+1})-\delta_{W}(u^{\kappa}_1).
\label{delta_U_delta_u}
\end{equation}

Set $\tilde\kappa=\kappa+\lfloor n^{\tilde\beta}\rfloor-1$. Let $Z_{\kappa}(x)$ be the number of vertices in $V_0^{\kappa}$ having at most $|W|p-x\sqrt{|W|p(1-p)\ln\tilde\kappa}$ neighbors in $W$. Clearly, by~(\ref{moivre_laplace}) and~(\ref{phi_approximation}), ${\sf E}Z_{\kappa}(x)=\tilde\kappa P(x)$, ${\sf Var}Z_{\kappa}(x)=\tilde\kappa P(x)(1-P(x))<{\sf E}Z_{\kappa}(x)$,
$$
 P(x)=(1+o(1))\int_{-\infty}^{-x\sqrt{\ln\tilde\kappa}}\frac{1}{\sqrt{2\pi}}e^{-t^2/2}dt=\frac{1}{x\sqrt{2\pi\ln\tilde\kappa}}\tilde\kappa^{-x^2/2}(1+o(1))=\frac{1}{\tilde\kappa}e^{\lambda}(1+o(1)),\text{ where}
$$
$$
 \lambda=\ln\frac{\tilde\kappa^{1-x^2/2}}{x\sqrt{2\pi\ln\tilde\kappa}}.
$$




First, let $\lambda=-\sqrt[4]{\ln n}$. Then, by a direct computation, we get $x=\sqrt{2}+\frac{\sqrt[4]{\ln n}}{\sqrt{2}\ln\tilde\kappa}(1+o(1))$ and ${\sf P}(Z_{\kappa}(x)\geq 1)\leq {\sf E}Z_{\kappa}(x)=e^{\lambda}=e^{-\sqrt[4]{\ln n}}$. 
Therefore,
\begin{multline*}
 {\sf P}(\exists\kappa\in\{1,\ldots,|V|-\lfloor n^{\tilde\beta}\rfloor+1\}\quad Z_{\kappa}(x)\geq 1)\leq \\
 {\sf E}Z_{1}(x)+\sum_{\kappa=2}^{|V|-\lfloor n^{\tilde\beta}\rfloor+1}\frac{1}{\kappa+\lfloor n^{\tilde\beta}\rfloor-1} e^{\lambda}(1+o(1)) =O\left(\ln n e^{-\sqrt[4]{\ln n}}\right).
\end{multline*}
Then, a.a.s. for every $\kappa\in\{1,\ldots,|V|-\lfloor n^{\tilde\beta}\rfloor+1\}$, 
\begin{equation}
\delta_W(u^{\kappa}_1)>|W|p-\sqrt{2p(1-p)|W|\ln\tilde\kappa}-\frac{\sqrt{|W|p(1-p)\sqrt{\ln n}}}{\sqrt{2\ln\tilde\kappa}}(1+o(1)).
\label{delta_1}
\end{equation}

Second, let $\lambda=\sqrt[4]{\ln n}$. Then $x=\sqrt{2}-\frac{\sqrt[4]{\ln n}}{\sqrt{2}\ln\tilde\kappa}(1+o(1))$ and ${\sf E}Z_{\kappa}(x)=e^{\sqrt[4]{\ln n}}$. From Chernoff inequality, 
$$
 {\sf P}(Z_{\kappa}(x)\leq h)\leq e^{-\frac{1}{2}e^{\sqrt[4]{\ln n}}(1+o(1))}=o\left(\frac{1}{n}\right).
$$
Then, for every $\kappa\in\{1,\ldots,|V|-\lfloor n^{\tilde\beta}\rfloor+1\}$, with probability $1-o\left(\frac{1}{n}\right)$,
\begin{equation}
\delta_{W}(u^{\kappa}_{h+1})<|W|p-\sqrt{2p(1-p)|W|\ln\tilde\kappa}+\frac{\sqrt{|W|p(1-p)\sqrt{\ln n}}}{\sqrt{2\ln\tilde\kappa}}(1+o(1)).
\label{delta_2}
\end{equation}
Finally, from~(\ref{delta_U_delta_u}),~(\ref{delta_1}) and the union bound applied to~(\ref{delta_2}), we get that a.a.s., for every $\kappa$,
$$
\delta_{W}(\mathbf{u}_{\kappa-1})-\delta_{W}(\mathbf{u}_{\kappa})=O\left(\frac{\sqrt{n}}{\sqrt[4]{\ln n}}\right).
$$

This implies that the probability of the required event from the statement of Claim~\ref{claim:algorithm} tends to 1 but does not imply that it equals $1-o(1/n)$. To prove the latter statement, let us consider a partition $V=V^1\sqcup\ldots\sqcup V^{\lfloor\ln n\rfloor}$ such that $||V^i|-|V^j||\leq 1$ and apply the above arguments for each of the sets $V^i$ (clearly, they work well even when $|V|\geq n^{\beta-o(1)})$. Then, the probability that desired $v_1,\ldots,v_h$ do not appear in each of the sets $V^1,\ldots,V^{\lfloor\ln n\rfloor}$ equals $[o(1)]^{\lfloor\ln n\rfloor}=o(1/n)$ as needed.

\subsubsection{Proof of Claim~\ref{claim:two_sets}} 
\label{claim:two_sets_proof}

Let $r=\lceil 2/\varepsilon\rceil+1$. Let $V^1\sqcup\ldots\sqcup V^r$ be a partition of $V$ such that $||V^i|-|V^j||\leq 1$. Let us prove that, for every $i$, the probability of non-existence of the desired $w_1,\ldots,w_{\ell}$ in $V^i$ is $O(1/n^{\varepsilon})$. Then we immediately get the statement of Claim~\ref{claim:two_sets} due to the independency of sets $E(G(n,p)|_{V^i\cup U})\setminus E(G(n,p)|_U)$, $i\in[r]$: the probability that there are no $w_1,\ldots,w_{\ell}$ in every $V^i$ is $O(1/n^{\varepsilon r})=o(1/n^2)$.


Let $W$ be the number of $\ell$-sets of vertices $w_1,\ldots,w_{\ell}\in V^i$ having exactly 
$$
\eta:=\left\lfloor\left({\ell\choose 2}+\ell|U|\right)p+\varphi\right\rfloor
$$ 
edges between them or going to $U$. By~(\ref{moivre_laplace_equal}), for every set $\mathbf{w}=\{w_1,\ldots,w_\ell\}$, the probability of the event $B_{w_1,\ldots,w_{\ell}}$ that $w_1,\ldots,w_{\ell}$ have $\eta$ edges between them or going to $U$ equals $\Theta\left(\frac{1}{\sqrt{n}}\right)$. Therefore, 
$$
{\sf E}W=\Theta\left(\frac{|V|^{\ell}}{\sqrt{n}}\right)=\Omega(n^{\varepsilon}).
$$
Also, by~(\ref{moivre_laplace_equal}), ${\sf P}(\mathrm{Bin}(n,p)=x)=O\left(\frac{1}{\sqrt{n}}\right)$ uniformly over all $x\in\{0,1,\ldots,n\}$. Therefore,
\begin{multline*}
{\sf Var}W\leq
{\sf E}W+\sum_{\mathbf{w},\mathbf{\tilde w}}{\sf P}(B_{\mathbf{w}}\cap B_{\mathbf{\tilde w}})\leq\\
{\sf E}W+\sum_{\tau=1}^{\ell-1} |V|^{2\ell-\tau}\sum_a{\sf P}(w_1,\ldots,w_{\tau}\text{ have }a\text{ edges between them or going to }U)\times\\
\times\left[{\sf P}\left(\mathrm{Bin}\left[(\ell-\tau)|U|+{\ell\choose 2}-{\tau\choose 2},p\right]=\eta-a\right)\right]^2\\
={\sf E}W+\sum_{\tau=1}^{\ell-1}|V|^{2\ell-\tau}\Theta\left(\frac{1}{n}\right)={\sf E}W+O\left(\frac{|V|^{2\ell-1}}{n}\right),
\end{multline*}
where the summation in the first line is over all distinct $\ell$-tuples $\mathbf{w},\mathbf{\tilde w}$  of vertices from $V^i$ with a non-empty intersection; the vertices $w_1,\ldots,w_{\tau}$ in the second line are arbitrary vertices of $V^i$.

Therefore, by the Chebyshev's inequality, 
$$
{\sf P}(W=0)\leq\frac{\mathrm{Var}W}{({\sf E}W)^2}=\frac{1}{{\sf E}W}+O\left(\frac{1}{|V|}\right)=O\left(n^{-\varepsilon}\right).\quad\Box
$$



\section{Proof of Theorem~\ref{AK}}

Denote $m(k)=\sqrt{(n-k)n\ln{n\choose k}}$ for $k\geq\lfloor\varepsilon n\rfloor$ and $m(k)=k\sqrt{\ln {n\choose k}}$ for $k<\lfloor\varepsilon n\rfloor$.

\subsection{Proof of Theorem~\ref{AK}.(ii)}

First, let $k\in\{\lfloor\varepsilon n\rfloor,\ldots,n-1\}$, $Q=3\sqrt{p}$ and $\mu=Qm(k)$.

Let $U$ be a $k$-vertex subset of $[n]$. Then, by the Chernoff inequality, the number of edges $e_U$ in $G(n,p)$ having at least one vertex outside $U$ does not belong to the interval
$$
\mathcal{I}_k:=\left(p\left(k(n-k)+{{n-k}\choose 2}\right)-\frac{\mu}{2}, p\left(k(n-k)+{{n-k}\choose 2}\right)+\frac{\mu}{2}\right)
$$
with probability at most $2e^{-\frac{\mu^2}{8\left(pn(n-k)+\frac{\mu}{6}\right)}}=e^{-\frac{9}{8}\ln{n\choose k}(1+o(1))}$ since $m=\sqrt{n(n-k)\ln{n\choose k}}<(n-k)\sqrt{n\ln n}=o(n(n-k))$. The expected number of $k$-vertex sets $U$ such that $e_U\notin\mathcal{I}_k$ is at most
$$
 {n\choose k}2e^{-\frac{\mu^2}{8\left(pn(n-k)+\frac{\mu}{6}\right)}}=e^{-\frac{1}{8}\ln{n\choose k}(1+o(1))}.
$$
Therefore, the probability that there exist $k\geq\lfloor\varepsilon n\rfloor$ and a $k$-vertex subset $U$ of $[n]$ such that $e_U\notin\mathcal{I}_k$ is at most 
$$
\sum_{k=\lfloor\varepsilon n\rfloor}^{n-1}e^{-\frac{1}{8}\ln{n\choose k}(1+o(1))}\leq \sum_{\ell=1}^8 n^{-\ell/8+o(1)}+n(1-\varepsilon)n^{-9/8+o(1)}\to 0\text{ as }n\to\infty.
$$
Since having an interval for the number of edges in induced graphs is equivalent to having an interval for the number of edges outside these induced graphs, we get that a.a.s., for every $k$ in the interval, the set of sizes of $k$-vertex induced subgraphs does not contain a full interval of length at least $\mu$, as desired.\\

Second, let $k\in\left\{\left\lceil\frac{1}{p}\ln n\right\rceil,\ldots,\lfloor\varepsilon n\rfloor-1\right\}$, $Q=3\sqrt{p}$ and $\mu=Qm(k)$ as well.

Let $U$ be a $k$-vertex subset of $[n]$. Then, by the Chernoff inequality, the number of edges $\tilde e_U$ in the induced subgraph $G(n,p)|_U$ does not belong to the interval
$$
\mathcal{J}_k:=\left(p{k\choose 2}-\frac{\mu}{2}, p{k\choose 2}+\frac{\mu}{2}\right)
$$
with probability at most 
$$
2e^{-\frac{\mu^2}{8\left(\frac{pk^2}{2}+\frac{\mu}{6}\right)}}\leq 2e^{-\frac{\mu^2}{8pk^2}}=e^{-\frac{9}{8}\ln{n\choose k}(1+o(1))}
$$ 
since $\frac{pk^2}{2}\geq\frac{\sqrt{p}k\sqrt{k\ln n}}{2}>\frac{\mu}{6}$. The expected number of $k$-vertex sets $U$ such that $\tilde e_U\notin\mathcal{J}_k$ is at most $e^{-\frac{1}{8}\ln{n\choose k}(1+o(1))}$. Therefore, the probability that there exist $k\geq\lfloor\varepsilon n\rfloor$ and a $k$-vertex subset $U$ of $[n]$ such that $\tilde e_U\notin\mathcal{J}_k$ is at most 
$$
\sum_{k=\left\lceil\ln n/p\right\rceil}^{\lfloor\varepsilon n\rfloor-1}e^{-\frac{1}{8}\ln{n\choose k}(1+o(1))}=e^{-\frac{1}{8}\ln{n\choose k}(1+o(1))}\to 0\text{ as }n\to\infty.
$$

Finally, for $k\in\left\{1,\ldots,\left\lceil\frac{1}{p}\ln n\right\rceil-1\right\}$ set $Q=\frac{1}{\sqrt{p}}$. Then, the number of edges of a $k$-vertex graph should belong to the interval $\{0,1,\ldots,{k\choose 2}\}$ of the length smaller than $\frac{k^2}{2}\leq k\sqrt{k\frac{1}{p}\ln \frac{n}{k}}<\frac{1}{\sqrt{p}}k\sqrt{\ln{n\choose k}}$ for $n$ large enough. The latter expression equals $Qm(k)$, and this finishes the proof.

\subsection{Proof of Theorem~\ref{AK}.(i)}

Note that, for $\varepsilon>0$ small enough, the case $k<\varepsilon n$ was already considered in~\cite{AlonKostochka}. Fix such an $\varepsilon<\frac{1}{4}$. 

Here, we consider three cases separately: 1) $k<n-n^{1/4}$,  2) $n-n^{1/4}\leq k\leq n-2$ and 3) $k=n-1$.

\subsubsection{$\varepsilon n\leq k<n-n^{1/4}$}

Let $q=\frac{\varepsilon\sqrt{\varepsilon p(1-p)}}{33}$.\\

Divide the set $\{1,\ldots,n-k+14\}$ into three `almost equal' parts $V_1,V_2,V_3$ (such that $||V_i|-|V_j||\leq 1$ for $i,j\in\{1,2,3\}$). Set $\tilde n=\tilde n(k)=k-14$ and let $V^*_{\tilde n}=\{n-\tilde n +1,\ldots,n\}$. Let $G_{\tilde n}$ be the induced subgraph of $G(n,p)$ on $V^*_{\tilde n}$.\\

We start with two technical statements.\\

\begin{claim}
A.a.s., for every integer $k\in[\varepsilon n,n-n^{1/4})$, in $G_{\tilde n}$ there are more than $\varepsilon(n-k)$ vertices having degrees greater than $\tilde n(k) p+\sqrt{\frac{1}{2}\tilde n(k)p(1-p)\ln(n/(n-k))}.$
\label{claim:many_large_degrees}
\end{claim}

{\it Proof.} Fix $k$ and let $Y_k$ be the number of vertices in $G_{\tilde n}$ having degrees greater than $\tilde n(k) p+\sqrt{\frac{1}{2}\tilde n(k)p(1-p)\ln(n/(n-k))}$. Set $\alpha=\frac{1}{2}\ln(n/(n-k))$ and apply Claim~\ref{claim:general_degrees}. Since 
$$
g_{\varepsilon}=(1-\varepsilon)\frac{n}{\sqrt{\pi\ln(n/(n-k))}}\left(\frac{n-k}{n}\right)^{1/4}>(1-\varepsilon)(n-k)\frac{[n/(n-k)]^{3/4}}{\sqrt{\pi\ln(n/(n-k))}}>\varepsilon(n-k),
$$
we get that ${\sf P}(Y_k<\varepsilon(n-k))=o(1/n^2)$.
From the union bound, Claim~\ref{claim:many_large_degrees} follows. $\Box$\\

\begin{claim}
A.a.s., for every integer $k\in[\varepsilon n,n-n^{1/4})$, in $G_{\tilde n}$ there are no vertices having degrees at least $\tilde n p+\sqrt{6\tilde np(1-p)\ln\tilde n}$.
\label{claim:no_large_degrees}
\end{claim}

{\it Proof.}  Fix $k$ and let $Z_k$ be the number of vertices in $G_{\tilde n}$ having degrees at least $\tilde n p+\sqrt{6\tilde np(1-p)\ln\tilde n}$. Then, by~(\ref{moivre_laplace})~and~(\ref{phi_approximation}), for $n$ large enough,
$$
{\sf E}Z_k\sim\tilde n\frac{1}{\sqrt{12\pi\ln\tilde n}}e^{-3\ln\tilde n}<\frac{1}{\tilde n^2}.
$$
Then, the desired property holds with probability at least $1-\sum\limits_{k\in[\varepsilon n,n-n^{1/4})}\frac{1}{k^2}=1-O(1/n)$. $\Box$\\


Let us describe an algorithm of finding $\tau\in\mathbb{N}$ and constructing sequences of subsets $U_1\subset\ldots\subset U_{\tau}$ in $V^*_{\tilde n}$ and $\tilde U_1\subset\ldots\subset\tilde U_{\tau}$ in $V_1$ such that (below, we denote $V^*_{\tilde n}[i]:=V^*_{\tilde n}\cup\tilde U_i\setminus U_i$, $i\in[\tau]$) a.a.s.
\begin{equation}
 e\left(G(\tilde n,p)|_{V^*_{\tilde n}[\tau]}\right)\leq e\left(G(\tilde n,p)\right)-qm(k)
\label{after_the_last_step}
\end{equation}
and, for every $i\in[\tau]$, $|U_i|=|U_{\tilde i}|=i$, 
\begin{equation} 
 -(\sqrt{2}+\sqrt{6})\sqrt{\tilde np(1-p)\ln\tilde  n}<e\left(G(\tilde n,p)|_{V^*_{\tilde n}[i]}\right)-e\left(G(\tilde n,p)|_{V^*_{\tilde n}[i-1]}\right)<
  -\sqrt{\frac{1}{2}\tilde n p(1-p)\ln\frac{n}{n-k}},
\label{one_step}
\end{equation}
where $U_0=\tilde U_0=\varnothing$.

It would mean that, up to an $(\sqrt{2}+\sqrt{6})\sqrt{\tilde np(1-p)\ln \tilde n}$-error, every value from 
\begin{equation}
\left(e(G(\tilde n,p))-qm(k),e(G(\tilde n,p))\right)
\label{big_interval}
\end{equation} 
is admissible by the number of edges in an induced $\tilde n$-vertex subgraph of $G(n,p)$.

Note that, having sequences of sets $U_1\subset U_2\subset\ldots$ and $\tilde U_1\subset\tilde U_2\subset\ldots$ with $|U_i|=|\tilde U_i|=i$ satisfying~(\ref{one_step}), the inequality~(\ref{after_the_last_step}) becomes true once 
\begin{equation}
\tau\geq\frac{qm(k)}{\sqrt{\frac{1}{2}\tilde n p(1-p)\ln(n/(n-k))}}.
\label{number_of_steps}
\end{equation}

Below, we show that our algorithm runs at least $\frac{\varepsilon(n-k)}{15}$ steps, and this immediately implies the inequality~(\ref{number_of_steps}): for $n$ large enough,
$$
 \frac{\varepsilon(n-k)}{15}>\frac{\varepsilon\sqrt{\varepsilon p(1-p)/2}\sqrt{n (n-k)\ln\left[\left(\frac{n}{n-k}\right)^{n-k}\right]}}{16\sqrt{\frac{1}{2}\varepsilon n p(1-p)\ln\frac{n}{n-k}}}\geq\frac{qm(k)}{\sqrt{\frac{1}{2}\tilde n p(1-p)\ln\frac{n}{n-k}}}.
$$

At step $1$, $U_1=\{v_1\}$,  where $v_1$ is a vertex having maximum degree in $G_{\tilde n}$. Consider the set $\mathcal{A}_1\subset V_1$ of vertices having at most $(\tilde n-1)p$  and at least
$$
R=(\tilde n-1)p-\sqrt{2\tilde np(1-p)\ln\tilde n}
$$
edges going to $V^*_{\tilde n}\setminus\{U_1\}$. Let $\tilde v_1\in\mathcal{A}_1$ (if $\mathcal{A}_1$ is non-empty; otherwise, the algorithm terminates), and $\tilde U_1=\{\tilde v_1\}$. 

Since a vertex from $V_1$ has at most $(\tilde n-1)p$ and at least $R$ neighbors in $V^*_{\tilde n}\setminus U_1$ with probability $1/2+o(1)$ (see Section~\ref{sec:pre}), the set $\mathcal{A}_1$ is non-empty with probability at least $1-(1/2+o(1))^{|V_1|}$.\\

Assume that, at step $1\leq i<\frac{\varepsilon(n-k)}{15}$, we construct the target sets $U_i,\tilde U_i$ having $i$ vertices. 

At step $i+1$, take a set $U_{i+1}=U_i\cup\{v_{i+1}\}$ of $i+1$ vertices having maximum degrees in $G_{\tilde n}$. Consider the set $\mathcal{A}_{i+1}\subset V_1\setminus \tilde U_i$ of vertices having at most $(\tilde n-1)p$ and at least $R$ edges going to $(V^*_{\tilde n}\cup \tilde U_i)\setminus U_{i+1}$. Let  $\tilde v_{i+1}\in\mathcal{A}_{i+1}$ (if $\mathcal{A}_{i+1}$ is non-empty; otherwise, the algorithm terminates), and $\tilde U_{i+1}=\tilde U_i\cup\{\tilde v_{i+1}\}$.

Let us prove that, with high probability, the set $\mathcal{A}_{i+1}$ is non-empty. Given an $(\tilde n-1)$-set, the probability that an outside vertex has at most $(\tilde n-1)p$ and at least $R$ neighbors in this set, equals $1/2+o(1)$ (see Section~\ref{sec:pre}). By the union bound, the probability that there exists an $i$-set $\tilde U$ in $V_1$ such that every vertex in $V_1\setminus \tilde U$ has either at least $(\tilde n-1)p$ or at most $R$ neighbors in $(V^*_{\tilde n}\cup\tilde U)\setminus U_{i+1}$ is at most
$$
 {|V_1|\choose i}(2+o(1))^{i-|V_1|}\leq 
 {|V_1|\choose \lfloor |V_1|/5\rfloor}(2+o(1))^{\lfloor |V_1|/5\rfloor-|V_1|}\leq \left(\frac{5e}{16}+o(1)\right)^{|V_1|/5}
$$
since $i<\frac{\varepsilon(n-k)}{15}<\frac{|V_1|}{5}$ (by definition, $|V_1|\geq\frac{n-k-16}{3}$ and $\varepsilon<\frac{1}{4}$).\\

Summing up, with probability at least $1-e^{-\Omega(n^{1/4})}$, for every $k\in[\varepsilon n,n-n^{1/4})$, the described algorithm works at least $\lceil\frac{\varepsilon(n-k)}{15}\rceil$ steps. By Claims~\ref{claim:many_large_degrees},~\ref{claim:no_large_degrees}, a.a.s.~for every $k$ in the range, it gives the desired sets.

Fix $i\in\{1,\ldots,\lceil\frac{\varepsilon(n-k)}{15}\rceil\}$ and consider the algorithm output $V^*_{\tilde n}[i]$. Notice that this set still has $k-14$ vertices.\\


From Claim~\ref{claim:algorithm} applied to $V=V_2$, $W=V^*_{\tilde n}[i]$, $\ell=4$ and $\beta=1/4$, it follows that a.a.s., there exists $a>0$ such that, for $k\in[\varepsilon n,n-n^{1/4})$ and any real $\alpha_0\in(3,7)$, we may find a set of $h=10$ vertices in $V_2$ having $ph\tilde n-\alpha_0\sqrt{\tilde n p(1-p)\ln\tilde n}+\xi$ neighbors in $V^*_{\tilde n}[i]$, where $|\xi|\leq a n^{1/2}(\ln n)^{-1/4}$, i.e. up to an $O(\sqrt{n}/\ln^{1/4}n)$-error, every value from~(\ref{big_interval}) is admissible (since $\sqrt{6}+\sqrt{2}<4$).\\

Finally, consider the set $V_3$. Let $\hat U$ be a union of $V^*_{\tilde n}[i]$ with a subset of $V_2$ having $10$ vertices.  It remains to prove that, a.a.s., for every $\gamma=\gamma(n)$  such that $0\leq\gamma\leq 4a n^{1/2}(\ln n)^{-1/4}$ and $(4(\tilde n+h)+6)p-\gamma$ is an integer, there exist vertices $w_1,\ldots,w_{4}$ in $V_3$ such that the number of edges in $G(n,p)|_{\hat U\cup\{w_1,\ldots,w_{4}\}}$ adjacent to at least one of $w_1,\ldots,w_{4}$ is exactly $(4(\tilde n+h)+6)p-\gamma$. But this follows immediately from Claim~\ref{claim:two_sets} applied to $V=V_3$, $U=\hat U$, $\ell=4$ and $\varphi=-\gamma$ and the union bound over all values of $\gamma$ and $k$.

\subsubsection{$n-n^{1/4}\leq k\leq n-2$}

The result immediately follows from the following three technical statements.

\begin{claim}
A.a.s., for every 
$$
 d\in I:=\biggl[(n-1)p-\sqrt{\frac{1}{5}n p(1-p)\ln n},(n-1)p+\sqrt{\frac{1}{5}n p(1-p)\ln n}\biggr],
$$
in $G(n,p)$, there are at least $n^{3/10}$ vertices having degree $d$.
\label{claim_degrees}
\end{claim}

{\it Proof}.  Fix $d\in I$ and let $Z$ be the number of vertices in $G(n,p)$ having degrees equal to $d$. Set $\tilde\alpha^2=\frac{1}{5}\ln n$ and apply Claim~\ref{claim:general_degrees}. We get that ${\sf P}\left(Z<\frac{1}{2\sqrt{2\pi p(1-p)}}n^{2/5}\right)=o(1/n^2)$. The desired statement follows from the union bound. $\Box$\\

For a subset $U\subset [n]$ of size $s$, let $\delta_0(U)=\delta(U)-({s\choose 2}+s(n-s))p$ be the difference between $\delta(U)$ and its expected value.\\

Let $q=\frac{1}{3}\sqrt{p(1-p)}$, $k\in[n-n^{1/4},n-2]$.\\

\begin{claim}
Assume that, in a graph $\mathcal{G}$ on $[n]$, for every $d\in I$, there are at least $n^{3/10}$ vertices having degree $d$.
Then, there exists a sequence $D_1\leq D_2\leq\ldots\leq D_{\kappa}$ such that $D_1>qm/2$, $D_{\kappa}<-qm/2$, for every $i\in[\kappa-1]$, $D_i-D_{i+1}\leq n^{1/4}$, and, in $\mathcal{G}$, there are sets of vertices $U_1,\ldots,U_{\kappa}$ of size $n-k-1$ having $\delta_0(U_i)=D_i$ for $i\in[\kappa]$.
\end{claim}

{\it Proof}. 
Since $n-k-1<n^{3/10}$, we can find $n-k-1$ vertices $v_1,\ldots,v_{n-k-1}$ having degrees equal to 
$$
d^*=\left\lfloor (n-1)p+\sqrt{\frac{1}{5}n p(1-p)\ln n}\right\rfloor.
$$
Clearly, for the set $U_1$ of these vertices and large enough $n$, the following holds:
$$
\delta_0(U_1)\geq(n-k-1)d^*-{n-k-1\choose 2}-\left[{n-k-1\choose 2}+(n-k-1)(k+1)\right]p>
$$
$$
 (p-1){n-k-1\choose 2}-(n-k-1)+(n-k-1)\sqrt{\frac{1}{5}n p(1-p)\ln n}>\frac{1}{6}(n-k)\sqrt{n p(1-p)\ln n}>\frac{qm}{2}.
$$

Let us shift degrees of every vertex of $U_1$ (starting from $v_1$) by a series of replacements and show that the desired sets $U_{2},\ldots,U_{\kappa}$ can be obtained on the way.

In the first sequence of steps, at every step, replace $v_1$ with a vertex having degree deg$(v_1)-1$. At the $i$th step, the set $U_i$ is produced. Once 
$$
\mathrm{deg}(v_1)=d_*:=\left\lceil(n-1)p-\sqrt{\frac{1}{5}n p(1-p)\ln n}\right\rceil,
$$
we proceed with replacing $v_2$ in the same way. In the final sequence of steps, we replace the vertex $v_{n-k-1}$. Once 
$\mathrm{deg}(v_{n-k-1})=d_*$, we stop and get a set $U_{\kappa}$ having 
$$
\delta_0(U_{\kappa})\leq (n-k-1)d_*-\left[{n-k-1\choose 2}+(n-k-1)(k+1)\right]p<
$$
$$
 p{n-k-1\choose 2}+(n-k-1)-(n-k-1)\sqrt{\frac{1}{5}n p(1-p)\ln n}<-\frac{1}{6}(n-k)\sqrt{n p(1-p)\ln n}<-\frac{qm}{2}.
$$
Clearly, for every $i\in[\kappa-1]$, $\delta_0(U_{i-1})-\delta_0(U_i)\leq n-k-2<n^{1/4}$. $\quad\Box$\\

Therefore, a.a.s., for every integer $k\in[n-n^{1/4},n-2]$ and every integer $D\in[-qm/2,qm/2]$, we may find an $(n-k-1)$-set $U\subset[n]$ such that $\delta_0(U)\in[D,D+n^{1/4})$. Therefore, we may occupy by the numbers of edges of $(k+1)$-sets every $n^{1/4}$-th number of an interval of length $qm$. It remains to prove that we may remove one vertex from such sets to get the whole interval, which is stated in the following claim.

\begin{claim}
A.a.s., for every integer $k\in[n-n^{1/4},n-2]$, every non-negative $d\leq n^{1/4}$ such that $pk+d$ is integer and every $(n-k-1)$-set $U\subset [n]$, there exists a vertex $z\in [n]\setminus U$ having exactly $pk+d$ neighbors in $[n]\setminus U$.
\end{claim}

{\it Proof.} Fix an integer $k\in[n-n^{1/4},n-2]$ and a non-negative $d\leq n^{1/4}$ such that $pk+d$ is integer. Let $U\subset [n]$ be an $(n-k-1)$-set. Without loss of generality, assume that $[n]\setminus U=\{1,\ldots,k+1\}$. For $\ell\in\{1,\ldots,k+1\}$, let $\mathcal{A}_{\ell}=\mathcal{A}_{\ell}(d)$ be the event that the vertex $v_{\ell}$ has exactly $pk+d$ neighbors in $[n]\setminus U$. We should estimate $ {\sf P}(\overline{\mathcal{A}_1}\cap\ldots\cap\overline{\mathcal{A}_{k+1}})$.\\

Divide the set $\{1,\ldots,k+1\}$ into $K:=\left\lfloor\frac{k+1}{\lfloor n^{3/4}\ln^5 n\rfloor}\right\rfloor$ sets $W_1,\ldots,W_K$ of the same size $\lfloor n^{3/4}\ln^5 n\rfloor$ (up to a remainder of a size less than $\lfloor n^{3/4}\ln^5 n\rfloor$ --- we remove it and do not consider it any more). 

Fix $i\in\{1,\ldots,K\}$. Let $\mathcal{S}_i$ be the event that all the degrees of the subgraph induced on $W_i$ are inside
$$
J:=\left((|W_i|-1)p-\sqrt{n},(|W_i|-1)p+\sqrt{n}\right).
$$
By the Chernoff bound, for every $i\in\{1,\ldots,K\}$, 
$$
 {\sf P}\left(\overline{\mathcal{S}_i}\right)\leq|W_i|e^{-\frac{n^{1/4}}{2\ln^5 n}(1+o(1))}=e^{-\frac{n^{1/4}}{2\ln^5 n}(1+o(1))}.
$$

Without loss of generality, assume that $W_i=\{1,\ldots,w\}$,  $w=\lfloor n^{3/4}\ln^5 n\rfloor$. For every {\it possible} graph $\mathcal{G}$ on $W_i$ having all degrees inside $J$ (we denote $\Gamma_i$ the set of all possible graphs), let $\mathcal{B}[\mathcal{G}]=\{G(n,p)|_{W_i}=\mathcal{G}\}$. Clearly, for such $\mathcal{G}$,
$$
 {\sf P}\left(\left.\overline{\mathcal{A}_1}\cap\ldots\cap\overline{\mathcal{A}_w}\right|\mathcal{B}[\mathcal{G}]\right)=
 {\sf P}\left(\overline{\mathcal{A}_1[\mathcal{G}]}\right)\cdot\ldots\cdot{\sf P}\left(\overline{\mathcal{A}_{w}[\mathcal{G}]}\right),
$$
where $\mathcal{A}_{\ell}[\mathcal{G}]$ is the event that the number of neighbors of $v_{\ell}$ in $[n]\setminus (U\cup W_i)$ equals $pk+d-\mathrm{deg}_{\mathcal{G}}(v_{\ell})$. By~(\ref{moivre_laplace_equal}), for some constant $c>0$, ${\sf P}\left(\left.\mathcal{A}_{\ell}[\mathcal{G}]\right|\mathcal{W}\right)\geq\frac{c}{\sqrt{n}}$. 

Finally, we get
\begin{multline*}
 {\sf P}\left(\overline{\mathcal{A}_1}\cap\ldots\cap\overline{\mathcal{A}_{k+1}}\right)\leq \\
 {\sf P}\left(\overline{\mathcal{A}_1}\cap\ldots\cap\overline{\mathcal{A}_{k+1}}\cap\{\exists i\in\{1,\ldots,K\}\,\,\mathcal{S}_i\}\right)+{\sf P}\left(\overline{\{\exists i\in\{1,\ldots,K\}\,\,\mathcal{S}_i\}}\right)\leq\\
 \sum_{i=1}^{K}{\sf P}\left(\bigcap_{\ell\in W_i}\overline{\mathcal{A}_\ell}\cap\mathcal{S}_i\right)+{\sf P}\left(\bigcap_{i=1}^K\overline{\mathcal{S}_i}\right)=
 \sum_{i=1}^{K}\sum_{\mathcal{G}\in\Gamma_i}{\sf P}\left(\left.\bigcap_{\ell\in W_i}\overline{\mathcal{A}_\ell}\right|\mathcal{B}[\mathcal{G}]\right){\sf P}\left(\mathcal{B}[\mathcal{G}]\right)+{\sf P}\left(\bigcap_{i=1}^{K}\overline{\mathcal{S}_i}\right)\leq\\
 \left(1-\frac{c}{\sqrt{n}}\right)^{w}\sum_{i=1}^{K}\sum_{\mathcal{G}\in\Gamma_i}{\sf P}\left(\mathcal{B}[\mathcal{G}]\right)+\mathrm{exp}\left(-\frac{n^{1/2}}{2\ln^{10}n}(1+o(1))\right)\leq\\
e^{-\frac{c w}{\sqrt{n}}}\sum_{i=1}^{K}{\sf P}\left(\mathcal{S}_i\right)+\mathrm{exp}\left(-\frac{n^{1/2}}{2\ln^{10}n}(1+o(1))\right)\leq\\
K \mathrm{exp}\left(-c n^{1/4}\ln^5 n(1+o(1))\right)=\mathrm{exp}\left(-c n^{1/4}\ln^5 n(1+o(1))\right).
\end{multline*}

Then, the probability that there exists an integer $k\in[n-n^{1/4},n-2]$, a non-negative $d\leq n^{1/4}$ such that $pk+d$ is integer and an $(n-k-1)$-set $U\subset[n]$ such that every vertex $z\in[n]\setminus U$ does not have exactly $pk+d$ neighbors in $[n]\setminus U$ is at most
$$
 \sqrt{n}n^{n^{1/4}}e^{-c n^{1/4}\ln^5 n(1+o(1))}\to 0\text{ as }n\to\infty.\quad\Box
$$

\subsubsection{$k=n-1$}

Let $q=\frac{\sqrt{p(1-p)}}{2}$. We should prove that a.a.s.~the set of sizes of $(n-1)$-vertex subgraphs of $G(n,p)$ contains a full interval of length $q\sqrt{n\ln n}$, or, equivalently, the set of degrees of $G(n,p)$ contains a full interval of the same size. But this follows from Claim~\ref{claim_degrees}.

\section{Acknowledgements}

The first author's research is partially supported by NSF Grants DMS-1500121 and DMS-1764123, Arnold O.~Beckman Research Award (UIUC Campus Research Board RB 18132) and
the Langan Scholar Fund (UIUC). The second authors's research is supported by the Ministry of Science and Higher Education of the Russian Federation in the framework of MegaGrant no 075-15-2019-1926.

We thank the referees for the careful reading and useful suggestions.

\end{document}